\documentclass[12pt]{article}

\usepackage{mathrsfs}
\usepackage[english]{babel}
\usepackage[latin1]{inputenc}
\usepackage{times}
\usepackage[T1]{fontenc}
\usepackage[normalem] {ulem}
\usepackage{color}
\usepackage[all]{xy}
\usepackage{graphicx}
\usepackage{amssymb}
\usepackage{epstopdf}
\usepackage{amsmath}
\usepackage{amsthm}

\newtheorem{lemma}{Lemma}[section]
\newtheorem{theorem}[lemma]{Theorem}
\newtheorem{proposition}[lemma]{Proposition}

\theoremstyle{definition}
\newtheorem{example}[lemma]{Example}

\newtheorem{remark}[lemma]{Remark}

\newcommand{\N}{\mathbb{N}}
\newcommand{\A}{\mathrm{Aut}}

\newcommand{\Sym}{\mathrm{Sym}}

\newcommand{\sL}{\mathscr{L}}

\newcommand{\sO}{\mathscr{O}}

\newcommand{\sQ}{\mathscr{Q}}
\newcommand{\sR}{\mathscr{R}}
\newcommand{\bI}{\mathbb{I}}
\newcommand{\bJ}{\mathbb{J}}
\newcommand{\bK}{\mathbb{K}}

\begin{document}

\renewcommand{\title}[2][]{
  \begin{center}
    \begin{Large}
      \begin{bf}
        #2
      \end{bf}
    \end{Large}
  \end{center}}

\renewcommand{\author}[2][]{
 \begin{center}
    \begin{sc}
      #2    \end{sc}
  \end{center}}
\newcommand{\address}[1]{
  \begin{center}
        #1
  \end{center}}

\renewcommand{\maketitle}{}


\title{Lobe, Edge, and Arc Transitivity of Graphs of Connectivity 1}

\author{Jack E.\,Graver and Mark E.\,Watkins}
\address{Department of Mathematics\\
Syracuse University\\
Syracuse, NY}

\begin{abstract}
We give necessary and sufficient conditions for lobe-transitivity of locally finite and locally countable graphs whose connectivity equals 1.  We show further that, given any biconnected graph $\Lambda$ and a ``code'' assigned to each orbit of $\A(\Lambda)$, there exists a unique lobe-transitive graph $\Gamma$ of connectivity 1 whose lobes are copies of $\Lambda$ and is consistent with the given code at every vertex of $\Gamma$.  These results lead to necessary and sufficient conditions for a graph of connectivity $1$ to be edge-transitive and to be arc-transitive.  Countable graphs of connectivity 1 the action of whose automorphism groups is, respectively, vertex-transitive, primitive, regular, Cayley, and Frobenius had been previously characterized in the literature.
\end{abstract}

\section{Introduction}

Throughout this article, $\Gamma$ denotes a connected simple graph.  Consider the equivalence relation $\cong$ on the edge-set $E\Gamma$ of $\Gamma$ whereby $e_1\cong e_2$ if the edges $e_1$ and $e_2$ lie on a common cycle of $\Gamma$.  A {\em lobe}\footnote{Following O.\,Ore \cite{Or62}, a {\em lobe} is a subgraph graph that either consists of a cut-edge with its two incident vertices or is a maximal biconnected subgraph.    We eschew the word ``block'' for this purpose, as it leads to ambiguity when discussing imprimitivity.}  is a subgraph of $\Gamma$ induced by an equivalence class with respect to $\cong$, and we let $\sL(\Gamma)$ denotes the set of lobes of $\Gamma$.  A vertex is a {\em cut vertex} if it belongs to at least two different lobes.  Connected graphs (other than $K_2$) have connectivity 1 if and only if they have a cut vertex.

Graphs of connectivity 1 whose automorphism groups have certain given properties have been characterized.  Those whose automorphism groups are, respectively, transitive, primitive, and regular were characterized in \cite{JuWa77}.   In particular, primitive planar graphs of connectivity 1 were characterized in \cite{WaGr04}\footnote{For a much shorter proof of non-existence of 1-ended planar graphs with primitive automorphism group, see \cite{SiWa12}.}.  Cayley graphs of connectivity 1 were characterized in \cite{Wa78}.  Graphs of connectivity 1 with Frobenius automorphism groups were characterized in \cite{Wa18}.  In the present work, we complete this investigation; we characterize graphs of connectivity 1 whose automorphism groups act transitively on their set of lobes.  As a consequence, we obtain characterizations of edge-transitive graphs and arc-transitive graphs of connectivity 1.  

The conditions for a graph of connectivity 1 to be lobe-transitive or to be vertex-transitive are independent;  such a graph may have either property or neither one or both.  Such is not the case for edge- and arc-transitivity.  In Section \ref{LobeTrans} we give necessary and sufficient conditions for a graph to be lobe-transitive.  We further show that, given any biconnected graph $\Lambda$ and a ``list'' of orbit-multiplicities of copies of $\A(\Lambda)$, one can construct a lobe-transitive graph of connectivity 1 all of whose lobes are isomorphic to $\Lambda$ and locally respects the given list.  We give necessary and sufficient conditions for a countable graph of connectivity 1 to be  edge-transitive in Section \ref{EdgeTrans} and to be arc-transitive in Section \ref{ArcTrans}.  As the sets of conditions for these latter two properties are more intertwined with lobe-transitivity than the characterization of vertex-transitivity (for graphs of connectivity 1), scattered throughout are examples that illustrate some algebraic distinctions among these various properties.

\section{Preliminaries}

Throughout this article, the symbol $\N$ denotes the set of positive integers.  The symbols $\bI,\bJ,$ and $\bK$ denote subsets of $\N$ of the form $\{1,2,\ldots,n\}$ or the set $\N$ itself; they appear as sets of indices.  All graphs (and their valences) in this article are finite or countably infinite.  The symbol $\delta_{i,j}$ (the so-called ``Kronecker delta'') assumes the value 1 if $i=j$ and $0$ if $i\neq j$.     For a graph $\Lambda$ and any subgroup $H\leq\A(\Lambda)$, the set of orbits of $H$ acting on $V\Lambda$ is denoted by~$\sO(H)$.

The set of lobes of a graph $\Gamma$ is denoted by $\sL(\Gamma)$. We let $\{\sL_k:k\in\bK\}$ denote the partition of $\sL(\Gamma)$ into isomorphism classes of lobes.   For given $k\in\bK$ and a lobe $\Lambda\in\sL_k$, we let $\sO(\A(\Lambda))=\{[V\Lambda]_j:j\in\bJ_k\}$, and we understand that if $\sigma:\Lambda\to\Theta$ is an isomorphism between lobes in $\sL_k$, then $\sigma([V\Lambda]_j)=[V\Theta]_j$ for all $j\in \bJ_k$.  Finally, for each $k\in\bK$ and $j\in \bJ_k$, we define the function $\tau_j^{(k)}:V\Gamma\to\N$ by 
\begin{equation}\label{oldtau}
\tau_j^{(k)}(v)=|\{\Lambda\in\sL_k:v\in[V\Lambda]_j\}|.
\end{equation}
For $\Lambda_0\in\sL(\Gamma)$ and $n\in\N$, we recursively define
$$\Gamma_0(\Lambda_0)=\Lambda_0,$$
$$\Gamma_{n+1}(\Lambda_0)=\bigcup\{\Lambda\in\sL(\Gamma):V\Lambda\cap V\Gamma_n(\Lambda_0)\neq\emptyset\}.$$

\begin{lemma}\label{extension} {\rm( \cite{JuWa77} Lemma 3.1)}  Let $\Lambda, \Theta\in\sL(\Gamma)$ and let $n\in\N$.  If for each $k\in\bK$ and $j\in\bJ_k$, the function $\tau_j^{(k)}$ is constant on $V\Gamma$, then any isomorphism $\sigma_n:\Gamma_n(\Lambda)\to\Gamma_n(\Theta)$ admits an extension to an automorphism $\sigma\in\A(\Gamma)$.
\end{lemma} 

This lemma, which will be invoked below, had be been used in \cite{JuWa77} to prove the following characterization of vertex-transitive graphs of connectivity 1. 

\begin{theorem}\label{VertTrans}{\em (\cite{JuWa77} Theorem 3.2)}  Let $\Gamma$ be a graph of connectivity 1.  A necessary and sufficient condition for $\Gamma$ to be vertex-transitive is that all the functions $\tau_j^{(k)}$ be constant on $V\Gamma$.
\end{theorem}

{\bf Notation.}  When all the lobes of the graph $\Gamma$ are pairwise isomorphic, that is, the index set $\bK$ has but one element, then in Equation (\ref{oldtau}) the index $k$ is suppressed; we simply replace $\bJ_k$ by $\bJ$ and $\tau_j^{(k)}$ by $\tau_j$.

\section{Lobe-transitivity}\label{LobeTrans}

Let $\Gamma$ be a graph of connectivity exactly 1.  It is immediate from the above definitions that  the edge-sets of the lobes of $\Gamma$ are blocks of imprimitivity of the group $\A(\Gamma)$ acting on $E\Gamma$.  Hence any automorphism of $\Gamma$ must map lobes onto lobes, and therefore, if $\A(\Gamma)$ is to act transitively on ${\mathscr{L}}(\Gamma)$, then all the lobes of $\Gamma$ must be pairwise isomorphic.  However, pairwise-isomorphism of the lobes alone is not sufficient for lobe-transitivity, even when every vertex of $\Gamma$ lies in the same number of lobes. 

Let us first dispense with trees; the proof is elementary and hence omitted.

\begin{proposition}  A finite or countable tree is lobe-transitive (and simultaneously, edge-transitive) if and only if there exist $n_1, n_2\in \N\cup\{\aleph_0\}$ such that every edge has one incident vertex of valence $n_1$ and the other of valence $n_2$.  If $n_1=n_2$, the tree is also arc-transitive.
\end{proposition}

For graphs of connectivity 1 other than trees, we have the following characterization of lobe-transitivity.

\begin{theorem}\label{Main} Let $\Gamma$ be a graph of connectivity 1, and let $\Lambda_0$ be an arbitrary lobe of $\Gamma$.  Let $\{ P_i:i\in\bI\}$ be the set of orbits of $\A(\Gamma),$ and let $\sQ=\{Q_j:j\in\bJ\}$ be the set of those orbits of the stabilizer in $\A(\Gamma)$ of $\Lambda_0$ that are contained in $\Lambda_0$.  Then necessary and sufficient conditions for the graph $\Gamma$ to be lobe-transitive are:
\begin{enumerate}
\item For each lobe $\Lambda\in\sL(\Gamma)$, there exists an isomorphism $\sigma_\Lambda: \Lambda_0\to\Lambda.$
\item For each $j\in\bJ$, there exists a function $\tau_j:V\Gamma\to\N\cup\{0,\aleph_0\}$ given by
\begin{equation}\label{tau}
\tau_j(v)=|\{\Lambda\in\sL(\Gamma): v\in\sigma_\Lambda(Q_j)\}|
\end{equation}
such that for all $i\in\bI$,
\begin{enumerate}
\item $\tau_j$ is constant on $P_i$, and
\item for all $v\in V\Gamma$, $\tau_j(v)>0$ if and only if $v\in\sigma_\Lambda(Q_j)\subset P_i$ for some lobe $\Lambda\in\sL(\Gamma)$.
\end{enumerate}
\end{enumerate}
\end{theorem}
 
{\em Proof.}   {\em (Necessity)}  Suppose that $\Gamma$ is lobe-transitive.  For each lobe $\Lambda\in\sL(\Gamma)$, there is an automorphism $\overline{\sigma}_\Lambda\in\A(\Gamma)$ that maps the fixed lobe $\Lambda_0$ onto $\Lambda$.  The restriction to $\Lambda_0$ of $\overline{\sigma}_\Lambda$ is an isomorphism $\sigma_\Lambda:\Lambda_0\to\Lambda$ that satisfies condition (1). 

For any lobe $\Lambda$, an automorphism $\alpha\in\A(\Gamma)$ is in the stabilizer of $\Lambda$ if and only if $\overline{\sigma}_\Lambda^{-1}\alpha\overline{\sigma}_\Lambda$ is in the stabilizer of $\Lambda_0$.  It follows that the partition $\{\sigma_\Lambda(Q_j):j\in\bJ\}$ of $V\Lambda$ is the set of orbits of the stabilizer of $\Lambda$ that are contained in $\Lambda$.  Furthermore, since the stabilizer of $\Lambda_0$ is a subgroup of $\A(\Gamma)$, the partition $\{\sigma_\Lambda(Q_j):j\in\bJ\}$ of $V\Lambda$ refines the partition $\{P_i\cap V\Lambda:i\in\bI\}$.  If for some indices $i$ and $j$, the vertex $v$ satisfies $v\in\sigma_\Lambda(Q_j)\subset P_i$, then for any lobe $\Theta$, the vertex $\sigma_\Theta\sigma_\Lambda^{-1}(v)$ lies in $P_i\cap\sigma_\Theta(Q_j)$.  This implies that, for all $j\in\bJ$, the function $\tau_j$ as given in Equation (\ref{tau}) is well-defined and satisfies condition 2(a).

Suppose that for an arbitrary index $i$, the vertex $v$ lies in $P_i$.  Since by Equation (\ref{tau}), $\tau_j(v)$ counts for each $j\in\bJ$ the number of lobes $\Lambda$ such that $v$ lies in $\sigma_\Lambda(Q_j)$, it follows that $\tau_j(v)$ is positive exactly when $\sigma_\Lambda^{-1}(v)\in Q_j\subset P_i$ holds.\smallskip

{(\em Sufficiency)}   Assume conditions (1) and (2).  To prove that $\Gamma$ is lobe-transitive, it suffices to prove that every isomorphism $\sigma_\Theta:\Lambda_0\to\Theta$ is extendable to an automorphism of $\Gamma$.  (Note that in this direction of the proof, $\sigma_\Theta$ is not presumed to be the restriction to $\Lambda_0$ of an automorphism $\overline{\sigma}_\Theta\in\A(\Gamma)$ but, in fact, it is.)

Fix a lobe $\Theta_0$ and a vertex $v\in V\Lambda_0$ and let $w=\sigma_{\Theta_0}(v)$.  For some $j\in\bJ$, the vertex $v$ lies in $Q_j$, and so $w\in \sigma_{\Theta_0}(Q_j)$.  Since both $\tau_j(v)$ and $\tau_j(w)$ are therefore positive,  both $v$ and $w$ lie in the same orbit $P_i$ of $\A(\Gamma)$ by condition 2(b).  By 2(a), $\tau_j$ is constant on $P_i$ for each $j\in\bJ$, and so there exists a bijection  $\beta_j$ from the set of lobes $\Lambda$ such that $v\in\sigma_\Lambda(Q_j)$ onto the set of lobes $\Theta$ such that $w\in\sigma_\Theta(Q_j)$.  Let $\Lambda_1$ be a lobe in the former set, and let $\Theta_1=\beta_j(\Lambda_1)$.   

Although  $v$ and $w$ lie in the images of the orbit $Q_j$ in lobes $\Lambda_1$ and $\Theta_1$, respectively, the vertices $\sigma_{\Lambda_1}^{-1}(v)$ and $\sigma_{\Theta_1}^{-1}(w)$ need not be the same vertex of $\Lambda_0$.  However, since both vertices lie in the same orbit $Q_j$ of $\Lambda_0$, there exists an automorphism $\alpha\in\A(\Lambda_0)$ such that $\alpha\sigma_{\Lambda_1}^{-1}(v)=\sigma_{\Theta_1}^{-1}(w)$.  Then $\sigma_{\Theta_1}\alpha\sigma_{\Lambda_1}^{-1}$ is an isomorphism from $\Lambda_1$ onto $\Theta_1$ that maps $v$ onto $w$ and therefore agrees with $\sigma_{\Theta_0}$ at the vertex $v$ common to $\Lambda_0$ and $\Lambda_1$. 

The amalgamation of  $\sigma_{\Theta_1}\alpha\sigma^{-1}_{\Lambda_1}$ with $\sigma_{\Theta_0}$ is an isomorphism from $\Lambda_0\cup\Lambda_1$ to $\Theta_0\cup\Theta_1$.  By repeating this same technique, we can extend $\sigma_{\Theta_0}$ to all lobes adjacent to $\Lambda_0$ and then to all of their adjacent lobes and inductively to all of $\Gamma$. $\hfill\qed$
\medskip

\begin{example}\label{LobeDelta} Suppose that the lobes of $\Gamma$ are copies of some biconnected, vertex-transitive graph and that every vertex of $\Gamma$ is incident with exactly $m$ lobes where $m\geq2$.  By Theorem \ref{VertTrans}, $\Gamma$ is  vertex-transitive.  By Theorem \ref{Main}, $\Gamma$ is lobe-transitive, with $\sO(\A(\Gamma))$ being the trivial partition (with just one big cell $P_1$).  Also $|\bJ|=1$ and $\tau_1(v)=m$ for all $v\in V\Gamma$.
\end{example}

\begin{remark}  There exists a ``degenerate'' family of lobe-transitive graphs $\Gamma$ of connectivity 1 that have but a single cut vertex.  For some cardinal $\mathfrak{K}\geq2$, consider a collection of $\mathfrak{K}$ copies of a biconnected graph $\Lambda_0$, and let $v_0\in V\Lambda_0$.  In the notation of Theorem \ref{Main}, let $\sigma_\Lambda(v_0)=v_\Lambda$ for each copy $\Lambda$ of $\Lambda_0$ in the collection.  We obtain $\Gamma$ by identifying $v_0$ and all the vertices $v_\Lambda$ and naming the new amalgamated vertex $w$, which forms a singleton orbit $\{w\}=P_1$ of $\A(\Gamma)$.  Clearly $\Gamma$ is lobe-transitive with connectivity 1.  Finally for all $j\in\bJ$, we have $\tau_j(w)=\delta_{1,j}\mathfrak{K}$, while $\tau_j(v)=1$ for all vertices $v\neq w$.   If $\Lambda_0$ has finite diameter, then $\Gamma$ has zero ends (see \cite{Ha73}) when $\Lambda_0$ is finite and $\mathfrak{K}$ ends when $\Lambda_0$ is infinite; if $\Lambda_0$ has infinite diameter, then $\Gamma$ has at least $\mathfrak{K}$ ends.  Other than the graphs just described, all countable lobe-transitive graphs of connectivity 1 are ``tree-like'' with $\aleph_0$ cut vertices and either 2 or $2^{\aleph_0}$ ends.
\end{remark}

\begin{theorem}\label{Build} Let $\Lambda_0$ be any biconnected graph.  Let $H\leq\A(\Lambda_0)$, let $\sQ=\sO(H)=\{Q_j:j\in\bJ\}$, and let $\sR=\{R_k:k\in\bK\}$ be a partition of $V\Lambda_0$ refined by $\sQ$. For each $k\in\bK$, define the function $\mu_k:\bJ\to \N\cup\{0,\aleph_0\}$ where $\mu_k(j)>0$ if and only if $Q_j\subseteq R_k$ and (to avoid the triviality of a single lobe) $\sum_{j\in\bJ}\mu_k(j)\geq2$ for at least one $k\in\bK$.  Then there exists (up to isomorphism) a unique lobe-transitive graph $\Gamma$ of connectivity 1 such that
\begin{enumerate}
\item for each lobe $\Lambda\in\sL(\Gamma)$, there exists an isomorphism $\sigma_\Lambda:\Lambda_0\to\Lambda$;
\item for each vertex $v\in V\Gamma$ and each $j\in\bJ$, we have
$$\mu_k(j)=|\{\Lambda\in\sL(\Gamma):\sigma_\Lambda^{-1}(v)\in Q_j\subseteq R_k\}|.$$
\end{enumerate}
\end{theorem}
{\em Proof.}  Let $\Lambda_0,\ H,\ \sQ$, and $\mu_k$ be as postulated.  Let $\Gamma_0=\Lambda_0$ from which we construct $\Gamma_1$ as follows.  

Let $v$ be any vertex of $\Lambda_0$.  For some $j,k$, it must hold that $v\in Q_j\subseteq R_k$, and so $\mu_k(j)>0$.  For each $\ell$ such that $Q_\ell\subseteq R_k$, we postulate the existence of $\mu_k(\ell)$ copies $\Lambda$ of $\Lambda_0$ (including $\Lambda_0$ itself when $\ell=j$) such that, if $\sigma_\Lambda: \Lambda_0\to \Lambda$ is an isomorphism, then some vertex in $\sigma_\Lambda(Q_\ell)$ is identified with the vertex $v$.  The graph $\Gamma_1$ is produced by repeating this process for each vertex of $\Lambda_0$.  
We again repeat this process starting at each vertex $w\in V\Gamma_1\setminus V\Gamma_0$, the only notational change being that, if specifically $w\in\sigma_\Lambda(Q_{j'})$ for some $j'\in\bJ$, then we consider the set $\sigma_\Lambda(Q_{j'})$ in $\Lambda$ (instead of $Q_j$ in $\Lambda_0$) to which $w$ belongs.  Thus we construct $\Gamma_2$.

Inductively, suppose that $\Gamma_n$ has been constructed for some $n\geq2$.  Let $w\in V\Gamma_n\setminus V\Gamma_{n-1}$, and so $w\in\sigma_\Lambda(Q_m)$ holds for some $m\in\bJ$ and a unique lobe $\Lambda\in\sL(\Gamma_n)\setminus\sL(\Gamma_{n-1})$.  Supposing that $Q_m\subseteq R_k$, we postulate the existence of $\mu_k(m)$ new copies of $\Lambda_0$ that share only the vertex $w$ with $\Gamma_n$ according to the above identification.  In this way we construct $\Gamma_{n+1}$.  Finally, let $\Gamma=\bigcup_{n=0}^\infty \Gamma_n$.  

It remains only to prove that $\Gamma$ so-constructed is lobe-transitive.  Let $\Theta$ be any lobe of $\Gamma$.  By the above construction, all lobes of $\Gamma$ are pairwise isomorphic, and so there exists an isomorphism $\sigma_\Theta:\Lambda_0\to\Theta$.  Starting with $\Gamma'_0=\Theta$ and by using the technique in the proof of Sufficiency in Theorem \ref{Main}, one constructs a sequence $\Gamma'_0,\Gamma'_1,\ldots$ so that for all $n\in\N$, we have $\Gamma_n'\cong\Gamma_n$, and $\sigma_\Theta$ is extendable to an isomorphism from $\Gamma_n$ to $\Gamma'_n$.  Thus $\Gamma\cong\bigcup_{n=0}^\infty\Gamma'_n$, and $\sigma_\Theta$ can be extended to an automorphism of $\Gamma$.$\hfill\qed$
\bigskip

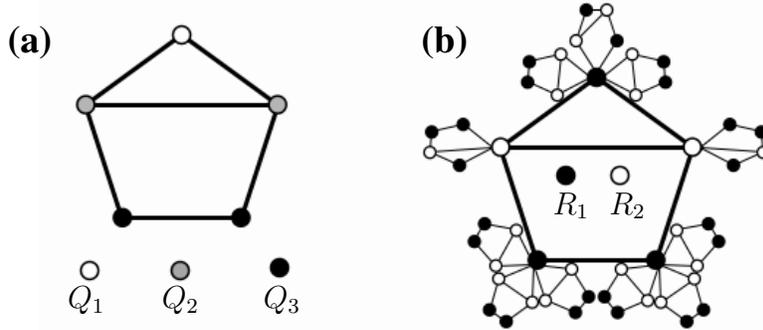
\begin{figure}[h]
\begin{picture}(15,3.4)(0,0)
\put(3.4,-.5){\includegraphics[width=10cm]{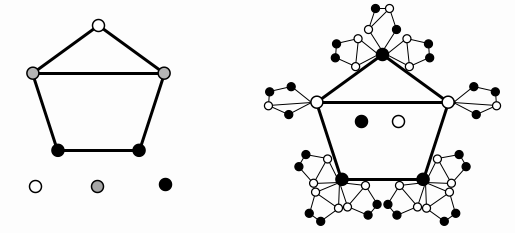}}
\put(6.4,-.1){$Q_3$}
\put(5.05,-.1){$Q_2$}
\put(3.8,-.1){$Q_1$}
\put(3,3.3){\large\bf (a)}
\put(8.5,3.3){\large\bf (b)}
\put(10.25,1.15){$R_1$}
\put(11,1.15){$R_2$}
\thicklines
\end{picture}
\caption{\label{5cyc}$\Lambda$ and $\Gamma_1$ from Example \ref{Chord5cyc}.}
\end{figure}

\begin{example}\label{Chord5cyc}  In the notation of Theorem \ref{Build}, let $\Lambda_0$ be the 5-cycle with one chord as shown in Figure \ref{5cyc}(a), and let $H=\A(\Lambda)$, yielding the orbit partition $\{Q_1,Q_2,Q_3\}$ as indicated.  Let $R_1=Q_1\cup Q_3$ and $R_2=Q_2$, giving $\bJ=\{1,2,3\}$ and $\bK=\{1,2\}$.  Define $\mu_1(1)=3,\ \mu_1(3)=1$, and $\mu_2(2)=2$.  Note that all other values of $\mu_1$ and $\mu_2$ must equal $0$.  Then $\Gamma_1$ is as seen in Figure~\ref{5cyc}(b).
\end{example}

The pairs of conditions in Theorems \ref{Main} and \ref{Build} may appear alike, but there is a notable difference between them.  This occurs when the arbitrarily chosen subgroup $H\leq\A(\Lambda_0)$ of Theorem \ref{Build} is a {\em proper} subgroup of the stabilizer of $\Lambda_0$ in $\A(\Gamma)$, where $\Gamma$ is the graph constructed from $\Lambda_0$ and the functions $\mu_k$ of Theorem \ref{Build}.   We illustrate this distinction with following example.
\begin{example}\label{distinctH}  Our initial lobe $\Lambda_0$ is a copy of $K_4$, with vertices labeled as in Figure \ref{ClotheslineBW2}(a), and so $\A(\Lambda_0)\cong\Sym(4)$ of order 24.  
We use $\Lambda_0$ to ``build'' the lobe-transitive graph $\Gamma$ shown four times in Figure \ref{ClotheslineBW2}(b).  The action on $\Lambda_0$ by the stabilizer of $\Lambda_0$ in $\A(\Gamma)$  is the 4-element group $\langle g_1\rangle\times\langle g_2\rangle$ whose generators have cycle representation $g_1=(v_1,v_2)$ and $g_2=(v_3,v_4)$.  The shadings of the vertices in the four depictions of $\Gamma$ in Figure \ref{ClotheslineBW2}(b) correspond respectively to the four different subgroups of $\langle g_1\rangle\times\langle g_2\rangle$ described below.  For the sake of simplicity, we assume $\sR=\sQ$.  

\begin{figure}[h]
\begin{picture}(16,6.5)(-.75,0) 
\put(0.2,-.5){\includegraphics[width=15cm]{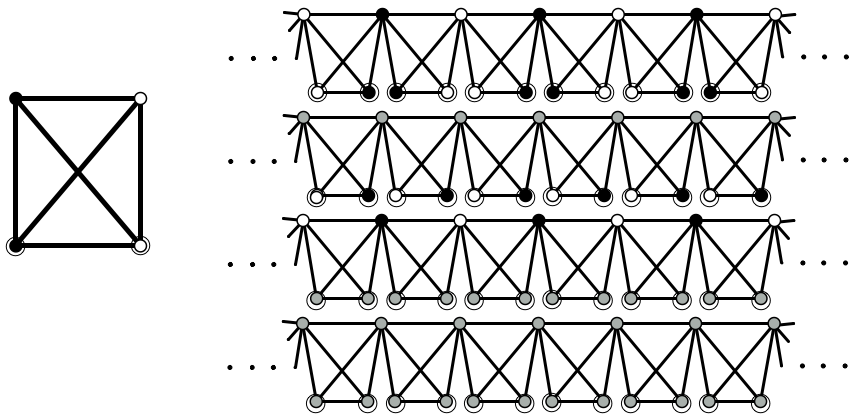}}
\put(-.1,5.3){$v_1$}
\put(2.85,5.3){$v_2$}
\put(-.1,2.8){$v_3$}
\put(2.85,2.8){$v_4$}
\put(0,6.4){\large \bf(a)}
\put(3.35,6.4){\large\bf (b)}
\put(4,6.1){(i)}
\put(4,4.25){(ii)}
\put(4,2.45){(iii)}
\put(4,0.6){(iv)}
\end{picture}
\caption{\label{ClotheslineBW2} The Clothesline Graph}
\end{figure}
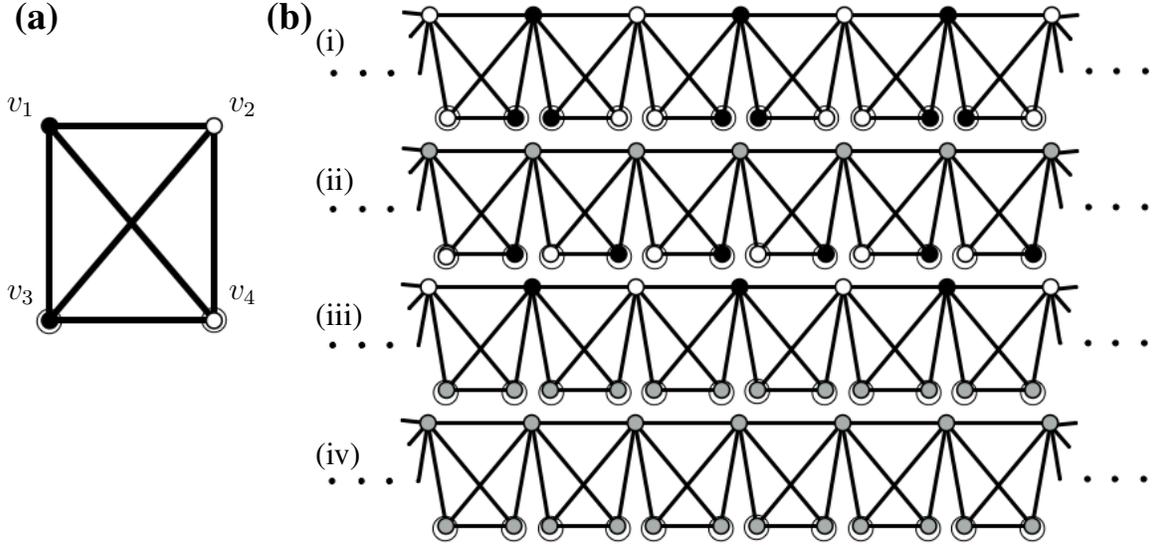

(i)  $H$ is the trivial group $\{\iota\}$.  Thus $H$ induces four orbits $Q_j=\{v_j\}$ for $j\in\{1,2,3,4\}$.  The functions $\mu_k$ are then given by $\mu_k(j)=2\delta_{j,k}$ for $k\in\{1,2\}$ and $\mu_k(j)=\delta_{j,k}$ for $k\in\{3,4\}$.

(ii) $H=\langle g_1\rangle$.    There are three orbits of $H:  Q_1=\{v_1,v_2\}, 
Q_2=\{v_3\}$, and $Q_3=\{v_4\}$, which give $\mu_1(1)=2, \mu_2(2)=\mu_3(3)=1$.  All other functional values are zero.

(iii)  $H=\langle g_2\rangle$.  Again there are three orbits of $H$ but not the same ones: $Q_1=\{v_1\}, Q_2=\{v_2\}$, and $Q_3=\{v_3,v_4\}$.  This gives $\mu_k(j)=2\delta_{j,k}$ for $k\in\{1,2\}$ and $\mu_3(j)=\delta_{j,3}$.

(iv) $H=\langle g_1\rangle\times\langle g_2\rangle$. Now there are just two orbits: $Q_1=\{v_1,v_2\}$ and $Q_2=\{v_3,v_4\}$.  Finally, $\mu_1(1)=2,\ \mu_2(2)=1$, and all other functional values are zero.

All four choices for $H$, the partition $\sQ$, and the functions $\mu_k$  clearly yield the same lobe-transitive graph $\Gamma$ of connectivity 1 by the construction of Theorem \ref{Build}
\end{example}

\section{Edge-transitivity}\label{EdgeTrans}

\begin{lemma}\label{EdgeLobe}  If $\Gamma$ is an edge-transitive (respectively, arc-transitive) graph, then $\Gamma$ is lobe-transitive and its lobes are also edge-transitive (respectively, arc-transitive). 
\end{lemma}
{\it Proof.}  For $i=1,2$, let $\Theta_i\in\sL(\Gamma)$, and let $e_i$ be an edge (respectively, arc) of $\Theta_i$.  There exists an an automorphism $\varphi\in\A(\Gamma)$ such that $\varphi(e_1)=e_2$.  Since $\varphi$ maps cycles through $e_1$ onto cycles through $e_2$, $\varphi$ must map $\Theta_1$ onto $\Theta_2$.  If $e_1$ and $e_2$ lie in the same lobe $\Theta$, then $\varphi$ leaves $\Theta$ invariant, and so its restriction to $\Theta$ is an automorphism of $\Theta$.
\hfill \qed 
\medskip

\begin{theorem}\label{Edges}  Let $\Gamma$ be a graph of connectivity 1 with more than one lobe, and let $\Lambda\in\sL(\Gamma)$.  Necessary and sufficient conditions for $\Gamma$ to be edge-transitive are the following:
\begin{enumerate}
\item The lobes of $\Gamma$ are edge-transitive. 
\item For each lobe $\Theta\in\sL(\Gamma)$, there exists an isomorphism $\sigma_\Theta:\Lambda\to\Theta$.
\item Exactly one of the following descriptions of $\Gamma$ holds:
\begin{enumerate}
\item Both $\Gamma$ and $\Lambda$ are vertex-transitive, in which case every vertex is incident with the same number $\geq2$ of lobes. 
\item The graph $\Gamma$ is vertex-transitive but $\Lambda$ is not vertex-transitive, in which case $\Lambda$ is bipartite with bipartition  $\{Q_1,Q_2\}$, and there exist constants $m_1,m_2\in\N\cup\{\aleph_0\}$ such that for $j=1,2$ and all $v\in V\Gamma$, it holds that $m_j=|\{\Theta\in\sL(\Gamma):v\in\sigma_\Theta(Q_j)\}|$.
\item The graph $\Gamma$ is not vertex-transitive, in which case $\Gamma$ is bipartite with bipartition $\{P_1,P_2\}$ and there exist constants $m_1,m_2\in\N\cup\{\aleph_0\}$, at least one of which is at least $2$, such that for $i=1,2$, if $v\in P_i$, then $|\{\Theta\in\sL(\Gamma):v\in\sigma_\Theta(P_j\cap V\Lambda)\}|=m_j\delta_{i,j}$.
\end{enumerate}
\end{enumerate}
\end{theorem}

{\em Proof.}  If the lobes of $\Gamma$ are edge-transitive and pairwise isomorphic and if $\Gamma$ satisfies any one of the three conditions, then it is immediate that $\Gamma$ is edge-transitive.

To prove necessity, suppose that $\Gamma$ is an edge-transitive graph of connectivity 1.  By Lemma \ref{EdgeLobe}, $\Gamma$ is also lobe-transitive and its lobes are edge-transitive, proving condition (1). Condition (2), which establishes notation for the remainder of this proof, also follows from Lemma \ref{EdgeLobe}.  

To prove (3), we continue the notation of Theorem \ref{Main} with $\bI$ being the index set for the set of orbits of $\A(\Gamma)$ and $\bJ$ being the index set for the orbits of the stabilizer of $\Lambda$ that are contained in $\Lambda$.  Since both $\Gamma$ and all of its lobes are edge-transitive, $|\bI|$ and $|\bJ|$ equal either 1 or 2.

If both $\Gamma$ and $\Lambda$ are vertex-transitive, then $|\bI|=|\bJ|=1$, and for every vertex $v\in V\Gamma$, $\tau_1(v)=m$ holds for some $m\geq2.$   This is case 3(a).

Since any odd cycle in $\Gamma$ would be contained in a lobe of $\Gamma$, it holds that $\Gamma$ is bipartite if and only if every lobe is bipartite.  If either $\Gamma$ or $\Lambda$ is not vertex-transitive, then each -- and hence both -- are bipartite, and the sides of the bipartitions (whether or not they are entire orbits of the appropriate automorphism group) are blocks of imprimitivity systems.  Let $\{P_1,P_2\}$ be the bipartition of $V\Gamma$, and so $\{P_1\cap V\Theta,P_2\cap V\Theta\}$ is the bipartition of any lobe $\Theta$.  Equivalently, letting $\{Q_1,Q_2\}$ denote the bipartition of $\Lambda$, we have $P_i=\bigcup_{\Theta\in\sL(\Gamma)}\sigma_\Theta(Q_i)$ for $i=1,2$.

Suppose that $\Gamma$ is vertex-transitive but $\Lambda$ is not, and so $|\bI|=1$ and $|\bJ|=2$.  By Theorem \ref{VertTrans}, there exist constants $m_1,m_2\in\N\cup\{\aleph_0\}$ such that for all $v\in V\Gamma$ and $j\in\bJ$, we have $m_j=\tau_j(v)=|\{\Theta\in\sL(\Gamma):v\in \sigma_\Theta(Q_j)\}|$.

Finally, suppose that $\Gamma$ is not vertex-transitive, and so $P_1$ and $P_2$ are the orbits of $\A(\Gamma)$.  Also, $\Lambda$ is bipartite with bipartition $\{Q_1,Q_2\}$, where $Q_i=P_i\cap V\Lambda$.  As no automorphism of $\Gamma$ swaps $P_1$ with $P_2$, no automorphism of $\Gamma$ swaps $Q_1$ with $Q_2$ (even though $\Lambda$ may be vertex-transitive!).  Hence $|\bI|=|\bJ|=2$.  It follows now from Theorem \ref{Main} that the functions $\tau_j$ satisfy the condition $\tau_j(v)>0$ if and only if $v\in P_j$.  That means that there exist constants $m_1,m_2\in\N\cup\{\aleph_0\}$, at least one of which is greater than 1, such that, if $v\in P_i$, then $\tau_j(v)=m_j\delta_{i,j}$.
\hfill\qed

\begin{example}\label{Kst}  Suppose in the notation of Theorem \ref{Edges} that $\Gamma$ is edge-transitive and $\Lambda$ is the complete bipartite graph $K_{s,t}$ with $|Q_1|=s$ and $|Q_2|=t$.  Suppose that every vertex of $\Gamma$ is incident with exactly two lobes.  If $s=t$, then both $\Gamma$ and $\Lambda$ are vertex-transitive, and we have case 3(a) of the theorem.  If $s\neq t$ and every vertex lies in one image of $Q_1$ and one image of $Q_2$, then we have the situation of case 3(b).  If again $s\neq t$ but each vertex lies in either two images of $Q_1$ or two images of $Q_2$, then we have the situation described in case 3(c). 
\end{example}

\begin{remark}  With regard to Example \ref{Kst}, we note that having $s=t$ does not assure vertex-transitivity of edge-transitive bipartite graphs. There exist edge-transitive, non-vertex-transitive, finite bipartite graphs where the two sides of the bipartition have the same size.  Such graphs are called {\em semisymmetric}.  The smallest such graph, on 20 vertices with valence 4, was found by J.\,Folkman \cite{Fo67}, who also found several infinite families of semisymmetric graphs.  Many more such families as well as forbidden values for $s$ were determined by A.\,V.\,Ivanov~\cite{Iv87}.
\end{remark}

\begin{example} This simple example illustrates how the converse of Lemma \ref{EdgeLobe} is false, even though the lobes themselves may be highly symmetric.   Let $\Gamma$ be a graph of connectivity 1 whose lobes are copies of the Petersen graph (which is 3-arc-transitive!).  For each lobe $\Lambda$,  let $[V\Lambda]_1$ and $[V\Lambda]_2$ denote the vertex sets of disjoint 5-cycles indexed ``consistently,'' i.e., if $\Lambda$ and $\Theta$ share a vertex $v$, then $v\in[V\Lambda]_i\cap[V\Theta]_i$ for $i=1$ or $i=2$.  (Observe that $\Gamma$ is not bipartite.)  For $i=1,2$ define  $P_i=\bigcup_{\Theta\in\sL(\Gamma)}[V\Theta]_i$, and suppose that  each vertex in $P_i$ belongs to exactly $m_i$ lobes.  The graph $\Gamma$ is lobe-transitive by Theorem \ref{Main} and is both vertex- and edge-transitive when $m_1=m_2$, but $\Gamma$ is neither vertex- nor edge-transitive when $m_1\neq m_2$.
\end{example}

\section{Arc-transitivity}\label{ArcTrans}  

 \begin{theorem}\label{Arcs}  
 Let $\Gamma$ be a graph of connectivity 1.  Necessary and sufficient conditions for $\Gamma$ to be arc-transitive are the following:
\begin{enumerate}
\item The lobes of $\Gamma$ are arc-transitive. 
\item The lobes of $\Gamma$ are pairwise isomorphic.
\item  All vertices of $\Gamma$ are incident with the same number of lobes.
\end{enumerate}
\end{theorem} 
{\em Proof.  (Necessity)}  Suppose that $\Gamma$ is arc-transitive.  Conditions (1) and (2) follow from Lemma \ref{EdgeLobe}.  Since arc-transitivity implies vertex-transitivity, condition (3) holds. 

{\em (Sufficiency)}  Assume that the three conditions hold.  For $m=1,2$, let $a_m$ be an arc of $\Gamma$, and let $\Theta_m$ be the lobe containing $a_m$.  By condition (2), there exists an isomorphism $\sigma:\Theta_1\to\Theta_2$.  By condition (1), there exists an automorphism $\alpha\in\A(\Theta_2)$ such that $\alpha\sigma(a_1)=a_2$.   By condition (3), the functions $\tau_j$ of Equation (\ref{oldtau}) are constant on $V\Gamma$.  (In fact, since the lobes are vertex-transitive, there is only one such function.)  It now follows from Lemma \ref{extension} that $\alpha\sigma$ is extendable to all of~$\Gamma$.\hfill\qed

\begin{remark}  If conditions (1) and (3) of Theorem \ref{ArcTrans} were replaced by {\em the lobes are edge-transitive and $\Gamma$ is vertex-transitive}, the {\em sufficiency} argument would fail.  There exist finite graphs \cite{Bo70} and countably infinite graphs with polynomial growth rate \cite{ThWa89} that are vertex- and edge-transitive but not arc-transitive.  Let $\Lambda$ denote such a graph, and consider a graph $\Gamma$ whose lobes are isomorphic to $\Lambda$ with the same number of lobes incident with every vertex.  Then $\Gamma$ itself is vertex- and edge-transitive but not arc-transitive.
\end{remark}
The following proposition is elementary.

\begin{proposition}  For all $k\geq2$, the only $k$-arc-transitive graphs of connectivity 1 are trees of constant valence.
\end{proposition}


\end{document}